\documentclass[12pt,a4paper,oneside]{article}
  \usepackage{authblk} %for writing authors
  \usepackage{url}
  \usepackage{amsthm}
  \usepackage{amsmath}
  \usepackage{amsfonts}
  \usepackage{amssymb}
  \usepackage{amscd}
  \usepackage{authblk}
  \usepackage{bm}
  \usepackage[mathscr]{eucal}%for \mathscr font
  \usepackage{mathtools} % for 'bsmallmatrix' environment
  \usepackage{physics}%for using Dirac's bra-ket notation, code: $\bra{\Psi}\ket{\Psi}$, $\expval{A}{\Psi}$, $\braket{x}{y}$
  \usepackage{color}
  \usepackage{fancyhdr}
  \usepackage{hyperref}
  \usepackage{marginnote}%for comments in margin, use \marginpar{text}
  \usepackage[top=2.5cm, bottom=2.7cm,left=2.5cm, right=2.5cm, marginparwidth=1.8cm]{geometry}
  \usepackage{framed}% use \begin{frame} text \end{frame} to put it in a framed box
  \usepackage[noabbrev]{cleveref} %for numbering in align* use \cref{...} for refering
  \usepackage[T1]{fontenc} %for writing multiple authors with multiple affiliations
 \usepackage[utf8]{inputenc}%for writing multiple authors
 \usepackage{blkarray}%for block matrices with subspaces marked for example see below

\usepackage{relsize} % example: \mathlarger{\mathlarger{\Pi}}
  \usepackage{enumerate}
   \usepackage{tikz}
 \usepackage{forest}
 \usepackage{qtree}  
 
  \numberwithin{equation}{section}
 
 \allowdisplaybreaks[1] 
  
  \theoremstyle{definition}  %Title and number in bold, body in normal font.
   \newtheorem{defn}{Definition}[section]

   \newtheorem{rmk}[defn]{Remark}

  \theoremstyle{plain}  %n Title and number in bold, body in italic (default).
   \newtheorem{thm}[defn]{Theorem}
   \newtheorem{lem}[defn]{Lemma}
   \newtheorem{prop}[defn]{Proposition}

  \theoremstyle{remark} %Title and number in italic, body in normal font.

   \newcommand{\BC}{\mathbb{C}}
   \newcommand{\BR}{\mathbb{R}}
   
   \newcommand{\BN}{\mathbb{N}}
   
   \newcommand{\ML}{{ ML^2(\BC\setminus\BR_-; q)}}

\newcommand{\llangle}{\langle \kern -0.2em \langle}	
\newcommand{\rrangle}{\rangle \kern -0.2em \rangle}

\usepackage{centernot}
\usepackage{todonotes}

 \newcommand{\numberthis}{\refstepcounter{equation}\tag{\theequation}}  
\makeatletter
\newcommand\footnoteref[1]{\protected@xdef\@thefnmark{\ref{#1}}\@footnotemark}
\makeatother

\let\tr\relax 
\DeclareMathOperator{\tr}{Tr}

\fontsize{12}{14}

\title{Thermal States on Mittag-Leffler Fock Space of the Slitted Plane}
\author[1]{Natanael Alpay}
\author[2]{Tiju Cherian John}

\affil[1]{Universiy of California Irvine, California, USA, nalpay@uci.edu}
\affil[2]{The University of Arizona, Tucson, Arizona, USA, tijucherian@fulbrightmail.org}
\date{}
\begin{document}
\maketitle
\begin{abstract}
Number states and thermal states form an important class of physical states in quantum theory.
A mathematical framework for studying these states is that of a Fock space over an appropriate Hilbert space. 
Several generalizations of the usual Bosonic Fock space have appeared recently due to their importance in many areas of mathematics and other scientific domains.
One of the most prominent generalization of Fock spaces is the Mittag-Leffler (ML) Fock space of the slitted plane. Natural generalizations of the basic operators of quantum theory can be obtained on  ML Fock spaces.
Following the construction of the creation and annihilation operators in the Mittag-Leffler Fock space of the slitted plane by 
Rosenfeld, Russo, and Dixon, (\href{https://doi.org/10.1016/j.jmaa.2018.03.036}{J. Math. Anal. Appl. 463, 2, 2018}).
We construct and study the number states and thermal states on the ML Fock space of the slitted plane. Thermal states on usual Fock space form an important subclass of the so called quantum gaussian states, an analogous theory of more general quantum states (like squeezed states and Bell states) on ML Fock spaces is an area open for further exploration.\\
\end{abstract}



\section{Introduction}

The Fock-Segal-Bargmann space or the Bosonic Fock space is a Hilbert space consisting 
of entire functions that are square integrable on the complex plane with respect to the normalized gaussian measure \cite{bargmann1961hilbert, Hall-98}.
%In 1961 Bargmann introduced a Hilbert space of entire functions 
The Fock space is the unique Hilbert space on which the creation and annihilation operators are adjoint of each other, i.e. $D^*=M_z$, where
$$\displaystyle M_zf(z):=zf(z)\quad \text{and} \quad Df(z):=\frac{d}{dz}f(z)$$ are closed, densely defined operators that satisfy the classical commutation rule 
$\left[D,M_z\right]=\mathcal{I}$,
where $\mathcal{I}$ is the identity operator. Recently, there were attempts in the literature to the generalize the construction of Fock space in different directions. For example,  the quaternionic Fock space and quaternionic Segal-Bargmann transform in the setting of slice hyperholomorphic functions were introduced in \cite{alpay2014fock,diki2017quaternionic}. The authors of \cite{alpay2022generalized}  studied a generalized Fock space built using different measures and found a new geometric description of the Fock space based on a special weight extending the classical normalized gaussian measure.

Another important construction of a generalized Fock space called Mittag-Leffler (ML) Fock space was introduced and studied by Rosenfeld et al. in \cite{mlf}.
They proved that classical complex monomials form an orthogonal basis of ML Fock space and computed their norms. The authors of \cite{mlf}, also defined analogues of the creation and annihilation operators on ML Fock spaces.
This new class of generalized Fock spaces led to many related works studying several classical objects in this more general setting. The associated Bargmann transform was studied in \cite{alpay2023mittag,williams2022segal}, fractional derivatives were studied in \cite{akgul2018novel,li2023hybrid}. In \cite{lian2021heisenberg} the authors established the Heisenberg’s uncertainty inequality associated with the Caputo fractional derivative. The associated  Littlewood-Paley identity  was studied in \cite{singh2021mittag}, some operator properties on the Mittag-Leffler spaces were studied in \cite{hai2021weighted}, and fractional oscillations were studied in \cite{li2022accurate}.
 In the present work, we take a quantum theoretic approach and study the analogues of number states and thermal states on the Mittag-Leffler space.

 Number states and thermal states are central objects in quantum optics, quantum information theory and more generally in quantum probability. Recall that by a state on a Hilbert space we mean a positive semidefinite trace class operator with unit trace.
Number states are vector states that form an orthonormal basis of the Fock space, meaning they are rank one projections corresponding to a sequence of orthornormal vectors.  
In the formalism of quantum optics, number states form a von Neumann measurement which counts the number of photons in a light pulse \cite{Sudarshan-Klauder}.
Thermal states are heavily used in  quantum optics and quantum information theory  \cite{Sudarshan-Klauder, weedbrook-et-al-2012}.
The relevance of number states and thermal states in general quantum mechanics and  continuous variable quantum information theory, can be seen in \cite{weedbrook-et-al-2012, Par12, ferraro05, WANG20071}.
For mathematically oriented expositions one may refer \cite{Par10, BhJoSr18, TCJ-thesis-2018, TCJ-KRP-2021}.

In this article, we generalize the number states and thermal states of the usual Fock space to the setting of the Mittag-Leffler Fock space of the slitted plane (see Remark \ref{rmk:number}, Definition \ref{defn:thermal} and Theorem \ref{thm:main}).
When restricted to the usual Fock space, our construction reduces to the number states and thermal states on the usual Fock space.
The theory of functions of a complex variable and special functions appear naturally in our constructions.
Several open problems also arise in this context.
\subsection{Discussion of the Main Result}
Now we provide an outline of our constructions in this article.
In the usual Bosonic Fock space, the number states are the eigenvectors of the number operator $a^{\dagger}a$ [where $a^{\dagger}$ denotes the creation operator, and $a$ denotes the annihilation operator]. 
The analogues of the number states in Mittag-Leffler Fock space is obtained by identifying an orthonormal sequence of eigenvectors of the operator $a^{\dagger}a$ in this space.
In the usual Fock space the thermal states with inverse temperature $s$ is defined as $$\rho_s=(1-e^{-s})e^{-sa^{\dagger}a},$$
where the scalar $(1-e^{-s})= \frac{1}{\operatorname{Tr} e^{-sa^{\dagger}a}}$. Following the definition of $a^{\dagger}$ and $a$ in the Mittag-Leffler Fock Space of the slitted plane, we prove that the positive semidefinite operator $e^{-sa^{\dagger}a}$ is a trace class operator on the Mittag-Leffler Fock space of the slitted plane, and thus we can define the thermal states of the ML space as 
$$ \frac{1}{\operatorname{Tr} e^{-sa^{\dagger}a}} e^{-s a^{\dagger}a }, $$ where $a^{\dagger}$ and $a$ above are respectively, the creation and annihilation operators of the ML Fock space of the slitted plane. Formal definitions and properties of the objects we require are provided in Section \ref{sec:prelim}.

Constructing the analogous of Weyl unitary operators and Bogoliubov Transformation of the usual Focks space are interesting open problems for future work.
More generally, the general structure of Mittag-Leffler spaces can be further explored if we can find analogous of Weyl  representation and metaplectic representation of the slitted plane into Mittag-Leffler Fock space

\section{Preliminaries}\label{sec:prelim}
In this section, we revise the important definitions and theorems that will be used throughout this article.
More specifically, we revise the notion of  ML space, which was first introduced in \cite{mlf}. 

\begin{defn}\label{defn:ml2cq}[ML space, 
{\cite[Definition 4 and Theorem 3.2]{mlf}}]
For $q>0$, the Mittag-Leffler space (sometimes referred to as the  Mittag-Leffler Fock space or MLF  space), $ML^2(\BC;q)$, of entire functions of order $q>0$ is defined in \cite{mlf} as 
 $$\displaystyle ML^2(\mathbb{C};q)=\left\lbrace  f(z)=\sum_{n=0}^{\infty}a_nz^n \;:\; \quad \frac{1}{q\pi}\int_{\mathbb{C}}|f(z)|^2|z|^{\frac{2}{q}-2}e^{-\frac{|z|^2}{q}}\dd z< \infty \right\rbrace,$$
 or equivalently using the sequential characterization
$$
M L^2(\mathbb{C} ; q)=\left\{f(z)=\sum_{n=0}^{\infty} a_n z^n\; :\; \sum_{n=0}^{\infty}| a_n|^2 \Gamma(q n+1)<\infty\right\}, .
$$ where $\Gamma(\cdot)$ is the usual Gamma function,    $\Gamma(z) := \int_0^\infty t^{z-1}e^{-t}$, for all $z\in \BC\setminus \{-1,-2,\dots\}$.
If $f(z) = \sum_{n=0}^{\infty} a_nz^{n}$ and $g(z)= \sum_{n=0}^{\infty} b_n z^{n}$, then \begin{equation}
    \label{scalar-product-ML2}
    \braket{f}{g}_{M L^2(\mathbb{C} ; q)}= \sum_{n=0}^{\infty} \bar{a}_nb_n\Gamma(qn+1).
\end{equation}
It is also the RKHS of entire functions associated with the kernel functions given by
$$
K_q(z, w)=E_q(\bar{z} w)=\sum_{n=0}^{\infty} \frac{\bar{z}^n w^n}{\Gamma(q n+1)},
$$
where $E_q(\cdot)$ is the Mittag-Leffler special function
\end{defn}

\begin{rmk}
    It can be noted, that for $q=1$, we obtain the classical geometric characterization Fock space $\mathcal{F}$ \cite{bargmann1961hilbert}. Also, following the convention in quantum theory, our scalar  product is linear in the second variable and anti-linear in the first variable
\end{rmk}

\begin{lem}[{\cite[Lemma 3.1]{mlf}}]\label{lm:31}
The set of functions  
$$
\left\{g_n(z)\right\}_{n=0}^{\infty}=\left\{\frac{z^n}{\sqrt{\Gamma(q n+1)}}\right\}_{n=0}^{\infty}
$$
is an orthonormal basis for $M L^2(\mathbb{C} ; q)$.
\end{lem}

\begin{defn}{\cite[Definition 7]{mlf}}\label{df:inner}
For $q>0$, the Mittag-Leffler space on the slitted plane is defined as
$$
M L^2\left(\mathbb{C} \backslash \mathbb{R}_{-} ; q\right):=\left\{f: \mathbb{C} \backslash \mathbb{R}_{-} \rightarrow \mathbb{C}: f\left(z^{1 / q}\right) \in M L^2(\mathbb{C} ; q)\right\}
$$
i.e.
$$
M L^2\left(\mathbb{C} \backslash \mathbb{R}_{-} ; q\right)=\left\{f: \mathbb{C} \backslash \mathbb{R}_{-}\rightarrow \BC, f(z)=\sum_{n=0}^{\infty} a_n z^{q n}: \sum_{n=0}^{\infty}\left|a_n\right|^2 \Gamma(q n+1)<\infty\right\}.
$$
The Mittag-Leffler space of the slitted plane is a RKHS equipped with the scalar product
$$
\braket{f(z)}{g(z)}_{M L^2(\mathbb{C} \backslash \mathbb{R}_{-} ; q)}=\braket{f\left(z^{1 / q}\right)} {g\left(z^{1 / q}\right)}_{M L^2(\mathbb{C} ; q)} .
$$ Hence if $f(z) = \sum_{n=0}^{\infty} a_nz^{qn}$ and $g(z)= \sum_{n=0}^{\infty} b_n z^{qn}$, then \begin{equation}
    \label{scalar-product-ML2-slitted-plan}
    \braket{f}{g}_{M L^2(\mathbb{C} \backslash \mathbb{R}_{-} ; q)}=\braket{f\left(z^{1 / q}\right)}{g\left(z^{1 / q}\right)}_{M L^2(\mathbb{C} ; q)} =\sum_{n=0}^{\infty} \bar{a}_nb_n\Gamma(qn+1).
\end{equation}
\end{defn}

 We recall the definition of Caputo fractional derivatives denoted by $D^q_*$ for $q>0$ in the ML space.
\begin{defn}[Caputo fractional derivative in ML spaces {\cite[ Definition 9]{mlf}}]\label{df:caputo}
Let $q>0$ and {$f\in ML^2(\mathbb{C} \setminus \mathbb{R};q)$ }
be given by the power series expansion $f(z) =\sum_{n=0}^{\infty}a_n z^{qn}$, then the Caputo fractional derivative of $f$ of order $q$ denoted by $D^{q}_{*}$ for $z\in \mathbb{C}\setminus \mathbb{R}_{-}$ is given by
\[ 
D^{q}_{*} f(z) = \sum_{n=1}^{\infty} a_n \frac{\Gamma(qn+1)}{\Gamma(q(n-1)+1)}z^{q(n-1)}.
\]
whenever 
\[ 
\sum_{n=1}^{\infty} |a_n|^2 \left|\frac{\Gamma(qn+1)}{\Gamma(q(n-1)+1)}\right|^2<\infty.
\]
\end{defn}

Finally, we recall the analogues of creation and annihilation  operators in the framework of ML space of the slitted plane

\begin{defn}[{\cite[Definition 10]{mlf}}]\label{def:aaddager}
    Define the following subspaces of $M L^2(\mathbb{C} \backslash \mathbb{R_{-}} ; q)$, \begin{align*}
\operatorname{Dom}\left(a^{\dagger}\right)&:=\left\{f \in M L^2\left(\mathbb{C} \backslash \mathbb{R}_{-};q\right): z^q f \in M L^2\left(\mathbb{C} \backslash \mathbb{R}_{-}; q\right)\right\}.\\
\operatorname{Dom}\left(a\right)&:=\left\{f \in M L^2\left(\mathbb{C} \backslash \mathbb{R}_{-};q\right): D_*^q f \in M L^2\left(\mathbb{C} \backslash \mathbb{R}_{-};q\right)\right\}
    \end{align*}
    Now the operators $a^{\dagger}$ and $a$ are defined on  as follows
    \begin{align*}
            \left(a^{\dagger} f\right)(z)&:=M_{z^q} f(z), \quad \forall f\in \operatorname{Dom}\left(a^{\dagger}\right)\\
             \left(a f\right)(z)&:=D_*^q f(z), \quad \forall f\in \operatorname{Dom}\left(a\right).
        \end{align*}
\end{defn}

\begin{prop}[{\cite[Proposition 5.2]{mlf}}]
    Let $q>0$. The operators $a^{\dagger}$ and $a$ are closed densely defined operators that are adjoint to each other, i.e., $(a^{\dagger})^*=a$, and $a^*=a^{\dagger}$.
\end{prop}

Recall that the elements $f$ of $ML^2(\BC\setminus\BR_-;q)$ are such that $f(z^{\frac{1}{q}})\in ML^2(\BC;q)$. Therefore, in the light of Lemma \ref{lm:31}, it is straight forward to see that the functions $\{g_n(z^q)\}_{n\geq 0}$ form an orthornomal basis of $ML^2(\BC\setminus\BR_-;q)$. Nevertheless, we provide a proof of this fact for completion.
\begin{prop}
    Let $\psi_n(z):=g_n(z^q)=\frac{z^{qn}}{\sqrt{\Gamma(qn+1)}}\in ML^2(\BC\setminus\BR_-;q)$ for $n=0,1,2\dots$ then the set $\{\psi_n\}_{n\geq 0}$ is an orthonormal basis (ONB) for the space $ML^2(\BC\setminus\BR_-; q)$.
\end{prop}
\begin{proof}  
    Notice first from Definition \ref{df:inner} that $\psi_n\in ML^2(\BC\setminus\BR_{-};q)$ for every $n\geq 0$.
    Now we show that the set is indeed orthonormal in $ML^2(\BC\setminus\BR_{-};q)$. 
    To this end, let $m\neq n$, then from Definition \ref{df:inner},  it follows that
    \begin{align*}
\braket{\psi_m}{\psi_n}_{ML^2(\BC\setminus\BR_{-};q)} &= \braket{\psi_m(z^{1/q})}{ \psi_n(z^{1/q})}_{ML^2(\BC\setminus\BR_{-};q)}\\
        &= \braket{\frac{z^m}{\sqrt{\Gamma(qm+1)}}}{\frac{z^n} {\sqrt{\Gamma(qn+1)}}}_{ML^2(\BC;q)}\\
        &=0,
    \end{align*}
    where the last line follows from Lemma \ref{lm:31}. Hence $\{\psi_n\}_{n\geq 0}$ is indeed orthogonal in the Hilbert space $ML^2(\BC\setminus\BR_{-};q)$. Now we will show that $\norm{\psi_n}_{ML^2(\BC\setminus\BR_{-};q)} =1$ for all $n\geq 0$. To this end,  note that
    \begin{align*}
        \left|\left|\frac{z^{qn}}{\sqrt{\Gamma(qn+1)}} \right|\right|^2_{ML^2(\BC \setminus\BR_{-};q)}
        &=\braket{\frac{z^{qn}}{\sqrt{\Gamma(qn+1)}}}{\frac{z^{qn}}{\sqrt{\Gamma(qn+1)}}}_{ML^2(\BC \setminus\BR_{-};q)}\\
        &=\braket{\frac{z^{n}}{\sqrt{\Gamma(qn+1)}}}{\frac{z^{n}}{\sqrt{\Gamma(qn+1)}}}_{ML^2(\BC;q)}\\
        &=\left|\left|\frac{z^{n}}{\sqrt{\Gamma(qn+1)}} \right|\right|^2_{ML^2(\BC;q)}\\
        &=1,
    \end{align*}
    by Lemma \ref{lm:31}.
    Hence $\{\psi_n\}_{n\geq0}$ is an orthonormal set. To show that it is an ONB of $ML^2(\BC\setminus\BR_{-};q)$, we need to show that $\{\psi_n\}_{n\geq0}$ spans $ML^2(\BC\setminus\BR_{-};q)$.
    Let $f\in ML^2(\BC\setminus\BR_{-};q)$. 
    By Definition \ref{df:inner}, $g(z)=f(z^{\frac{1}{q}})\in ML^2(\BC;q)$. 
    By Lemma \ref{lm:31}
    we get
    \begin{align*}
        g(z)=f(z^{\frac{1}{q}})&=\sum_{n=0}^{\infty}a_n g_n(z),\quad 
    \end{align*}
    where  $\sum_{n=0}^{\infty}|a_n|^2\Gamma(qn+1)<\infty$ and $g_n(z)=\frac{z^n}{\sqrt{\Gamma(qn+1)}}$ with $\{g_n(z)\}_{n\geq 0}$ is an ONB of $ML^2(\BC;q)$. Recall from Definition \ref{defn:ml2cq} of $ML^2(\BC;q)$,  that $g(z)=f(z^{\frac{1}{q}})$ is an entire function, so
    \begin{align*}
        f(z)&=g(z^q) =\sum_{n=0}^{\infty} a_n g_n(z^{q})\\ &=\sum_{n=0}^{\infty} a_n\underbrace{\frac{z^{qn}}{\sqrt{\Gamma(qn+1)}}}_{=\psi_n(z)}\\
        &=\sum_{n=0}^{\infty} a_n\psi_n(z), \quad \forall z\in \BC\setminus \BR_-.
    \end{align*}
    Since $\{\ket{\psi_n}\}_{n\geq 0}$ is orthonormal in $ML^2(\BC\setminus\BR_{-};q)$, we have  
    \[ f=\sum_{n=0}^{\infty}a_n \ket{\psi_n}, \]
   in the norm of $ML^2(\BC\setminus\BR_{-};q)$,  which completes the proof.
\end{proof}
\begin{rmk}\label{rmk:number} Following the notation $\ket{n}$ used to denote  $n$-particle Fock state (number state) in the usual Fock space; for $q>0$, we use the notation $\ket{n_q}$ to denote the function $\psi_n$ in the previous lemma, i.e.,
    \begin{equation}
    \label{eq:nq}
    \ket{n_q}:= \ket{\frac{z^{qn}}{\Gamma(qn+1)}},\quad n\in \BN.
    \end{equation}
    and for $n=0$, $\ket{0_q}\equiv 1$.
    In particular, for $q=1$, $\ket{n_1} = \ket{n}$ in the usual Fock space. Thus the states $\{\ket{n_q}\}$ are a natural generalization of usual number states in the setting of ML Fock spaces.
\end{rmk}
\begin{prop}\label{prop:a-dagger-a}
The  orthonormal basis $\{\ket {n_q}\}_{n=0}^\infty$ is contained in the domain of $a^\dagger a$ and \begin{equation}
    \label{eq:a-dagger-a-on-nq}
    a^{\dagger}a \ket{n_q} =
    \begin{cases}
        \frac{\Gamma(qn+1)}{\Gamma(q(n-1)+1)}\ket{n_q} &n\geq 1\\
        0 &n=0
    \end{cases}.
\end{equation}  In particular, the operator $a^{\dagger}a$ initially defined on the span of the orthonormal basis $\{\ket{n_q}\}$ extends as a positive selfadjoint operator, again denoted by $a^\dagger a$ on its maximal domain.  
Furthermore, if $f\in ML^2(\BC\setminus\BR_-; q)$ is such that $f(z)= \sum_{n=0}^{\infty}a_n z^{qn}$, and $f\in {\rm{Dom}}(a^\dagger a)$ then, \[\left(a^{\dagger}a f\right) (z)  = \sum_{n=1}^{\infty}a_n\frac{\Gamma(qn+1)}{\Gamma(q(n-1)+1)}z^{qn}.\]   
\end{prop}
\begin{proof}
    First we show that $\ket{n_q}$ is in the domain of $a$, 
    to this, note that by Definition \ref{def:aaddager} 
    \[ D(a) = \{f\in \ML\mid D^q_*(f)\in\ML\} \]
    We will compute $D^q_*\ket{n_q}$ and show that $D^q_*\ket{n_q}\in \ML$.
    By Definition \ref{df:caputo} for $z\in \mathbb{C}\setminus\mathbb{R}_{-}$ we have
       \begin{align*}
        D_*^q\ket{n_q}(z)
        &=\begin{cases}
            \ket{(n-1)_q}(z),&n\geq 1\\
            0,&n=0
        \end{cases}.
    \end{align*}
    Hence  \begin{align*}
        a\ket{n_q}
        &=\begin{cases}
            \ket{(n-1)_q},&n\geq 1\\
            0,&n=0
        \end{cases}\in\ML.
    \end{align*}
   Now we prove that $a\ket{n_q} \in \operatorname{Dom} (a^\dagger)$, to this end note that 
        \[ a\ket{n_q}(z) =\frac{z^{(n-1)q}}{\Gamma(q(n-1)+1)},\quad z\in \BC\setminus R_- ,  \]
        so \[z^qa\ket{n_q}(z) = \frac{z^{qn}}{\Gamma(q(n-1)+1)} \in \ML.\]
        Therefore, by Definition \ref{def:aaddager} $a \ket{n_q}$ is in the domain of $a^{\dagger}$.
    Now  \begin{align*}
        a^\dagger a\ket{n_q}(z) =z^q(a\ket{n_q}) (z)= 
\frac{z^{qn}}{\Gamma(q(n-1)+1)}=\frac{\Gamma(qn+1)}{\Gamma(q(n-1)+1)}\underbrace{\left(\frac{z^{qn}}{\sqrt{\Gamma(qn+1)}}\right)}_{\ket{n_q}}.
    \end{align*}
    Therefore, the expression for $a^\dagger a\ket{n_q}$ is given by
     \begin{align*}
    a^{\dagger}a\ket{n_q}
    =\frac{\Gamma(qn+1)}{\Gamma(q(n-1)+1)}{\ket{n_q}}, \quad n\in \BN\cup \{0\}
    \end{align*}
    proving \eqref{eq:a-dagger-a-on-nq}. Since $\{\ket{n_q}\}$ forms an orthonormal basis of eigenvectors for $a^{\dagger}a$, we see that the operator $a^{\dagger}a$ initially defined on this orthonormal basis using the eigenvalue relation above extends to a closed positive operator on a maximal domain. 
    
    To get the general form 
   for $ a^{\dagger}a f$ where $f(z)= \sum_{n=0}^{\infty}a_n z^{qn}\in {\rm{Dom}}(a^\dagger a)$, note that \[f(z) = \sum_{n=0}^{\infty}a_n z^{qn} = \sum_{n=0}^{\infty}a_n \Gamma(qn+1) \frac{z^{qn}}{\Gamma(qn+1)},  \] i.e., \[\ket{f} = \sum_{n=0}^{\infty}a_n \Gamma(qn+1)\ket{n_q}.\]  Now by the closedness of the operator $a^\dagger a$,
    \begin{align*}
        a^{\dagger}a\ket{f} 
        &= \sum_{n=0}^{\infty}a_n \Gamma(qn+1)a^{\dagger}a\ket{n_q}\\
        &= \sum_{n=1}^{\infty} a_n\frac{\Gamma(qn+1)^2}{\Gamma(q(n-1)+1)} \ket{n_q}.
    \end{align*}
    This last expression is same as the one given in the statement of the present proposition.
\end{proof}
\begin{rmk}
 By   Proposition \ref{prop:a-dagger-a}, we see that the operator $a^\dagger a$ on the ML-space of the slitted plane, $ML^2(\BC; q)$,  is a diagonal operator with respect to the orthonormal basis $\{\ket{n_q}\}$. 
\end{rmk}

\section{Thermal States on ML-space of the Slitted Plane}
Following the theory of thermal states on a classical Fock space \cite{Par10, Par12, TCJ-KRP-2021, TCJ-thesis-2018} we define the thermal states of the $ML^2(\BC\setminus \BR_-; q)$ space as follows.
\begin{defn}\label{defn:thermal}
    Let $s>0$, a thermal state with parameter $s$ of $ML^2(\BC\setminus\BR_-; q)$ is defined as the operator \[\rho_s=\frac{e^{-sa^{\dagger}a}}{\tr e^{-sa^{\dagger}a} }.\]
\end{defn} 
Note that our definition of thermal states assumes that the operator $e^{-sa^{\dagger}a}$ defined on $ML^2(\BC\setminus \BR_-; q)$ is trace class. We will devote the rest of this section to prove that this is indeed the case. Once we prove that  the operator $e^{-sa^\dagger a}$ defined on $ML^2(\BC\setminus \BR_-; q)$ is trace class, we get a one parameter generalization of the usual thermal states on Fock spaces (i.e., the case $q=1$).
\begin{rmk}\label{rmk:nq}
We use the notation
$$n_q := \frac{\Gamma(qn+1)}{\Gamma(q(n-1)+1)}, \quad n\in \BN.$$ It is a fact that $\Gamma(x)>0$, for $x>0$ \cite[Page 273]{Gorenflo-Kilbas-Mainardi}, hence $n_q>0$ for all $n$. With this notation and Proposition \ref{prop:a-dagger-a} we have 
    \[a^{\dagger} a\ket{n_q} = n_q \ket{n_q}, \quad \forall n \in \BN\cup \{0\}.\]
    Therefore, to prove that the operator $e^{-sa^{\dagger} a}$ is trace class, we need to show that the series\[\sum_{n=0}^\infty e^{-sn_q}<\infty.\]
    Indeed, we prove this and it is the main technical result of our article and we achieve this in Theorem \ref{thm:main}. Before that we discuss some tools from analysis we need to prove our main theorem.
\end{rmk}

\begin{rmk}\label{rmk:digamma} In what follows we need to use the Digamma function. We summarize the properties important to us in this remark.
The Digamma function is the logarithmic derivative of the Gamma function, given by $$\psi(z)=\frac{\Gamma^{\prime}(z)}{\Gamma(z)}=\frac{\mathrm{d}}{\mathrm{d} z} \ln \Gamma(z)$$
It is a strictly increasing and concave function on the interval $(0, \infty)$ \cite[Section 6.3]{Abramowitz-Stegun-1964}. Furthermore, the asymptotic behavior is given by 
$$\psi(z) \sim \ln z-\frac{1}{2 z}+ O (|z|^{-2})$$
for large values of $z$ in the sector $|\arg z|<\pi-\varepsilon$, where $\varepsilon$ is a small positive constant {\cite[Page 259, Equation 6.3.18]{Abramowitz-Stegun-1964}} .
\end{rmk}
\begin{lem}\label{nq-lim}  For $q>0$, define $f_q:[1,\infty)\rightarrow \BR$ by \[f_q(x)=x_q:=\frac{\Gamma(qx+1)}{\Gamma(q(x-1)+1)}, \quad x>1.\] Then we have \begin{enumerate}
    \item\label{item:1} $f_q(n) = n_q$, where $n_q$ is as in Remark \ref{rmk:nq};
    \item \label{item:2} $f_q$ is a positive increasing function on $[1,\infty)$;
    \item \label{item:3}$f_q$ satisfies \[\lim_{x\to  \infty}f_q(x) = \infty\]
\end{enumerate}
\end{lem}
\begin{proof} Note that \ref{item:1} is a trivial remark. 

\ref{item:2}. Recall  that the Gamma function is strictly positive and differentiable on the positive real axis (See \cite[Appendix A.1.2]{Gorenflo-Kilbas-Mainardi}). Hence  $f_q(x)>0$ for all $x\in [1,\infty)$. 
Now we show that the function $f_q$ is monotonically increasing by showing that the derivative of $f_q$ is positive. 
    Using the Digamma function $\psi (x) = \frac{\Gamma'(x)}{\Gamma(x)}$ defined in Remark \ref{rmk:digamma}, we have,

    \begin{align*}
        \dv{x}f_q(x) &= \dv{x}\left( \frac{\Gamma(qx+1)}{\Gamma(q(x-1)+1)} \right) \\
        &=\frac{q\Gamma'(qx+1)\Gamma(q(x-1)+1) -q\Gamma'(q(x-1)+1)\Gamma(qx+1)}{(\Gamma(q(x-1)+1))^2}\\
        &=q\frac{\frac{\Gamma'(qx+1)}{\Gamma(qx+1)}\Gamma(qx+1)\Gamma(q(x-1)+1) -\frac{\Gamma'(q(x-1)+1)}{\Gamma(q(x-1)+1)}\Gamma(q(x-1)+1)\Gamma(qx+1)}{(\Gamma(q(x-1)+1))^2}\\
        &= \frac{q\Gamma(qx+1)(\psi(qx+1)-\psi(q(x-1)+1))}{\Gamma(q(x-1)+1)}\\
        %&= \frac{q qn\Gamma(qn)(\psi(qn+1)-\psi(q(n-1)+1))}{q(n-1)\Gamma(q(n-1))}\\
        %&= q x_q \left(\psi(qx+1)-\psi(q(x-1)+1)\right)\\
        &= \underbrace{q x_q}_{>0}\cdot\underbrace{(\psi(qx+1)-\psi(qx+1-q))}_{>0} \numberthis \label{eq:derivative_nq}\\
        &>0,
    \end{align*}
    where we used the positivity of the Gamma function on $(0,\infty)$ and the fact that Digamma function is increasing on $(0,\infty)$ as stated in Remark \ref{rmk:digamma}.
    Hence we indeed have a strictly increasing sequence. 
    
    \ref{item:3}. Now we show that $x_q\to \infty$ as $x\to\infty$.
We have 
$$
\frac{\Gamma(z+\alpha)}{\Gamma(z+\beta)}=z^{\alpha-\beta}\left[1+\frac{(\alpha-\beta)(\alpha+\beta-1)}{2 z}+O\left(|z|^{-2}\right)\right],
$$   for arbitrary constants $\alpha$ and $\beta$, and $z\in \BC$ with $|\arg z| \leqslant \pi-\delta$ for a some positive $\delta$ (see \cite[Appendix A.1.8]{Gorenflo-Kilbas-Mainardi} or \cite[Chapter 1, Exercise 9]{Lebedev-1965}). Now the lemma follows by taking $z=qn, \alpha = q$ and $\beta = 0$.
\end{proof}
\begin{rmk}\label{abscissa} We need to use the abscissa of convergence of a general Dirichlet series in our next theorem. We provide the basics of this notion in this remark.
    The \textit{general Dirichlet series} is an infinite series of the form 
    \begin{equation*}
        \sum_{n=1}^{\infty} a_ne^{-\lambda_n s}
    \end{equation*}    
    where $a_n,s$ are complex numbers and $\{\lambda_n\}_{n\in\mathbb{N}}$ is a strictly increasing sequence of nonnegative real numbers that goes to infinity. The abscissa of convergence of a Dirichlet series is given by 
    $$
\sigma_c=\inf \left\{\sigma \in \mathbb{R}: \sum_{n=1}^{\infty} a_n e^{-\lambda_n s} \text { converges for every } s \text { for which } \operatorname{Re}(s)>\sigma\right\} .
$$ 
From {\cite[Theorem 8]{hardy-reiz-1915}},  $\sigma_c$ is given by
$$
\sigma_c=\limsup _{n \rightarrow \infty} \frac{\log( |a_1|+|a_2|+\cdots+|a_n|)}{\lambda_n}.
$$
Meaning that the Dirichlet series converges for every $s$ such that ${\rm{Re}}(s)>\sigma_c$.
\end{rmk}
Now we prove the main theorem of this article which shows that the operator  $e^{-sa^\dagger a}$ is a trace class operator on $ML^2(\BC\setminus\BR_-;q)$. Combining our next theorem with  Remark \ref{rmk:nq} completes the construction of thermal states on ML Fock space of the slitted plane. 
\begin{thm}\label{thm:main}
For $q>0$ and $s>0$, the series $\sum_{n=0}^{\infty}e^{-sn_q}$ converges, i.e., \begin{equation}
\label{eq:e_sum}
\sum_{n=0}^{\infty}e^{-sn_q}<\infty,
\end{equation} 
where $n_q = \frac{\Gamma(nq+1)}{\Gamma((n-1)q+1)}$.
\end{thm}
\begin{proof} By Lemma \ref{nq-lim} the sequence $\{n_q\}$ is increasing, so 
the series $\sum_{n=0}^{\infty}e^{-sn_q}$  has the form of a general Dirichlet series as in Remark \ref{abscissa}, with $a_n=1$, $\lambda_n =n_q,  {n\in\mathbb{N}}$.  
Therefore, the abscissa of convergence for this series is given by \begin{equation}\label{eq:abscissa}
    \sigma_c=\limsup _{n \rightarrow \infty} \frac{\log n }{n_q}
\end{equation} We will show that for $q>0$, $\sigma_c= 0$ so that the series is convergent for all $s>0$ .
    In fact, we will show that $\lim_{n\to\infty}\frac{\log{n}}{n_q}=0$. To this end,  note that  by Lemma \ref{nq-lim}, $n_q\to\infty$ as $n\to \infty$, and since $\log n\to \infty$ as $n\to \infty$,    we can apply l'H\^opital rule, and show instead that 
    $$\lim_{x\to\infty}\frac{\dv x \log{x}}{\dv x f_q(x)} =0,$$ where $f_q$ is as in Lemma \ref{nq-lim}.
    Now by using the formula for the derivative of $f_q$ given in equation \eqref{eq:derivative_nq}, we have 
    \begin{align}\label{eq:lopital}
\lim_{x\to\infty}\frac{{\dv x}\log{x}}{{\dv x}f_q(x)} 
        &= \lim_{x\to\infty}
        \frac{\frac{1}{x}}{q x_q \left(\psi(qx+1)-\psi(q(x-1)+1)\right)\nonumber}\\
        &= \lim_{x\to\infty}\frac{1}{q x_q x\left(\psi(qx+1)-\psi(q(x-1)+1)\right)}.
    \end{align}
        Using the asymptotic behavior of $\psi$ given in Remark \ref{rmk:digamma},   we have \begin{align*}
            &\lim_{x\to \infty} x\left(\psi(qx+1)-\psi(q(x-1)+1)) \right) \\ &\phantom{...........}= \lim_{x\to \infty}x\left(\log\left( {1-\frac{q}{qx+1}} \right) +{ \frac{1}{2(qx+1)} - \frac{1}{2(qx+1)-2q}}\right)\\
            &\phantom{...........}= -1
        \end{align*}
Therefore,   we  have the limit in \eqref{eq:lopital} equal to $0$ because $x_q\rightarrow \infty$ as $x\rightarrow \infty$ by Lemma~\ref{nq-lim}. Thus we have shown that the abscissa of convergence given by \eqref{eq:abscissa} is $0$ which completes the proof.
\end{proof}
\begin{rmk}
Theorem \ref{defn:thermal} states that the operator $e^{-sa^\dagger a}$ is a trace class operator for all $s>0$. Thus Definition \ref{defn:thermal} of thermal states makes sense.
\end{rmk}
\section{Conclusion and Discussion}
Following the construction of the creation and annihilation operators on the Mittag-Leffler (ML) Fock space of the slitted plane by Rosenfeld, Russo and Dixon \cite{mlf}, it is natural to look for an analogous quantum theory on  ML spaces. As described in the introduction of this article we can already see such attempts in the literature. In this article, we construct ML counter parts of the thermal states (Definition \ref{defn:thermal} and Theorem \ref{thm:main}) and number states (Remark \ref{rmk:number}) of quantum theory. The theory of thermal states and number states of usual Fock space and their generalizations are very important in quantum optics and quantum information theory. Therefore, further investigations and generalizations of our work is interesting. For example, one can investigate ML analogues of coherent states and ask whether there is a symmetry group  associated with the thermal states and coherent states of ML spaces \cite{Par13}. Furthermore, studying quantum stochastic calculus on ML spaces can reveal deeper understanding of the subject.

\subsection*{Acknowledgements:}
 TCJ thanks the Fulbright Scholar Program and United States-India Educational Foundation for providing funding and other support to conduct this research through a Fulbright-Nehru Postdoctoral Fellowship (Grant number: 2594/FNPDR/2020), he also acknowledges the United States Army Research Office MURI award on
Quantum Network Science, awarded under grant number W911NF2110325 for partially funding this research. Both the authors gratefully acknowledge the support provided by the the conference \textit{Advances in Operator Theory with Applications to Mathematical Physics}, November 14-18, 2022,
Chapman University, Orange, CA, where part of this work was conducted.

\bibliographystyle{IEEEtran}
\bibliography{bibliography}
\end{document}